\documentclass[final,sort&compress,3p,times]{elsarticle}
\usepackage[]{amsmath}
\usepackage{mathrsfs}
\usepackage{multirow}

\journal{}
\begin{document}
\begin{frontmatter}

\title{A direct algorithm of one-dimensional optimal system for the group invariant solutions}

\author[label1]{Xiaorui Hu\corref{cor1}}
\ead{lansexiaoer@163.com}
\author[label2]{Yuqi Li}
\author[label2]{Yong Chen\corref{cor1}}
\ead{ychen@sei.ecnu.edu.cn}

\address[label1]{Department of Applied Mathematics, Zhejiang University of Technology, Hangzhou 310023, China}
\address[label2]{Shanghai Key Laboratory of Trustworthy Computing,
East China Normal University, Shanghai 200062, China}

\cortext[cor1]{Corresponding author.}

\begin{abstract}

A direct and systematic algorithm is proposed to find one-dimensional optimal system for the group invariant solutions, which is attributed to the classification of its corresponding one-dimensional Lie algebra. Since the method is based on different values of all the invariants, the process itself can both guarantee the comprehensiveness and demonstrate the inequivalence of the optimal system, with no further proof. To illustrate our method  more clearly , we give a couple of well-known examples: the Korteweg-de Vries (KdV) equation and the heat equation.

\end{abstract}

\begin{keyword}
group invariant solutions;Lie algebra; optimal system; invariants; the adjoint  transformation  matrix
\end{keyword}

\end{frontmatter}

\section{Introduction}

Symmetry group theory built by Sophus Lie  plays an important role in constructing explicit solutions for integrable and non-integrable nonlinear equations. For any given subgroup, an original nonlinear system can be reduced to a system with fewer independent variables which corresponds to  group invariant solutions.
Since there are almost always an infinite amount of such subgroups, it is usually not feasible to list all possible group invariant solutions to the system. It is anticipated to find all those inequivalent group invariant solutions, that is to say, to give them a classification. The problem of classifying the subgroups and reduction to optimal systems takes on more importance for multidimensional PDEs. Given a  group that leaves a PDE invariant, one desires to minimize the search for
group-invariant solutions to that of finding inequivalent
branches of solutions, which leads to the concept of the
optimal systems. Consequently, the problem of determining the optimal system of subgroups is  reduced to the corresponding problem for  subalgebras.
In applications, one usually constructs the optimal system of subalgebras, from which the optimal systems of subgroup and group invariant solutions are reconstructed.

The adjoint representation of a Lie group on its Lie algebra was known to
Lie. Its use in classifying group-invariant solutions appears in Ovsiannikov~\cite{Ovsiannikov}. Ovsiannikov demonstrated
the construction of the one-dimensional optimal system for the Lie algebra, using a global matrix for the adjoint transformation and  sketched the construction of higher-dimensional optimal systems with a simple example.
The method has received extensive development by Patera, Winternitz and Zassenhaus~\cite{Patera1, Patera2} and many examples of optimal systems
of subgroups for the important Lie groups of mathematical physics were obtained. In the investigation of the connections between Lie group and special functions, Weisner~\cite{Weisner} firstly gave the classification of the symmetry algebra of the heat equation. For the higher-dimensional optimal systems of Lie algebra,  Galas~\cite{Galas} also developed Ovsiannikov's idea by  removing equivalent subalgebras and the problem of a nonsolvable algebra was also discussed, which is generally harder than that for a solvable algebra. Some examples of optimal systems can also be
found in Ibragimov~\cite{Ibragimov1,Ibragimov2}.

Here we are concerned with the one-dimensional  optimal  system of subalgebras. For the one-dimensional optimal systems, the technique of Ovsiannikov has been used until Olver gives a slightly  different and  elegant technique. Olver~\cite{Olver} constructed
a table of adjoint operators to simplify a general element  in Lie algebra as much as possible and applied the technique to the Korteweg-de Vries (KdV) equation and the heat equation. Since it only depends on fragments of the theory of Lie algebras, Olver's method as developed here has the feature of being very elementary. Based on Olver's method, we have also constructed many interesting and important invariant solutions~\cite{hu1,hu2,hu3,hu4,hu5}  for a number of systems of PDEs in atmosphere and geometric field.  However,  as Olver said, although some sophisticated techniques are available for Lie algebras with additional structure, in essence this problem is attacked by the na\"{i}ve approach. One knows that, if one calls a list of $\{\tilde{v}_{\alpha}\}_{\alpha \in \mathcal{A}}$ is an one-dimensional  optimal system, it must satisfy two conditions: (1) completeness--any one-dimensional subalgebra is equivalent to some $\tilde{v}_{\alpha}$;
(2) inequivalence--$\tilde{v}_{\alpha}$ and $\tilde{v}_{\beta}$ are inequivalent for distinct $\alpha$ and $\beta$. Despite these numerous results on finding the representatives of subalgebras, they did not illustrate that how these representatives are comprehensive and mutually inequivalent. Recently, Chou and Qu~\cite{Qu1,Qu2,Qu3} offer many numerical invariants to address the inequivalence among the elements in the optimal system.

The purpose of this paper is to give a systematic method for finding an optimal system of one-dimensional Lie algebra, which can both guarantee the comprehensiveness and the inequivalence. The idea is inspired by the observation that the killing form of the Lie algebra is an invariant  for the adjoint representation ~\cite{Olver}. Olver also points out that the detection of such an invariant is important since it places restrictions on how far one can expect to simplify the Lie algebra. In spite of the importance of the invariants for the Lie algebra, to the best of our knowledge, there are few literatures to use more common invariants except the killing form in the process of constructing optimal system. The purpose of this paper is to introduce a direct and valid method for providing  all the general invariants which are different from  the numerical invariants appearing in~\cite{Qu1,Qu2,Qu3} and then make the best use of them with the  adjoint matrix to classify subalgebras.  We shall demonstrate the new technique by
treating a couple of illustrative examples.

This paper is arranged as follows.
In section \ref{algorithm}, a direct algorithm of one-dimensional optimal system for the general symmetry algebra is proposed. Since the realization of our new algorithm
builds on  different  invariants and the adjoint matrix, a valid method for computing all the invariants is also given in this section.
In section \ref{heat equation}, we apply the new algorithm to a couple of well-known examples, i.e. the Korteweg-de Vries (KdV) equation and the heat equation, and construct their one-dimensional optimal systems step by step.
Finally, a brief conclusion is given in section \ref{conclusion}.

\section{A direct algorithm of one-dimensional optimal system}\label{algorithm}
Consider the $n$-dimensional symmetry algebra $\mathcal{G}$ of a differential
system, which is generated by the vector fields $\{v_1,v_2, \cdots, v_n\}$. The corresponding symmetry group of $\mathcal{G}$ is denoted as $G$.
Following Ovsiannikov~\cite{Ovsiannikov}, one calls two elements $v=\sum\limits_{i=1} ^{n} a_iv_i$ and $w=\sum\limits_{j=1} ^{n} b_j v_j$ in $\mathcal{G}$ equivalent if they satisfy one of the following conditions:

(1) one can find some transformation $g \in G$ so that
$Ad_g(w)=v$;

(2) there is $v=cw$ with $c$ being constant.

Here $Ad_g$ is the adjoint representation of $g$ and $Ad_g(w)=g^{-1}wg$. It needs to  note that the second condition is less obvious in  all the references but here it will play an important role in our method. The main tools used in our algorithm
are all the invariants and the adjoint matrix.

\subsection{Calculation  of  the  invariants} \label{invariants}

A real function $\phi$ on  the Lie algebra  $\mathcal{G}$ is called an invariant if $\phi(Ad_g(v))=\phi(v)$ for all $v \in \mathcal{G}$ and all $g \in G$. If two vectors $v$ and $w$ are equivalent under the adjoint action, it is necessary that $\phi(v)=\phi(w)$ for any invariant $\phi$. If we let $v=\sum\limits_{i=1} ^{n} a_iv_i$, then the invariant $\phi$ can be regarded as a function of $a_1, a_2, \ldots, a_n$.  As Olver said, the detection of such an invariant is important since it places restrictions on how far we can expect to simplify $v$. However, it is a pity that people did not care more invariants except the killing form. Now we will propose a valid method to find all the invariants of symmetry algebra and further make the best use of them to construct one-dimensional optimal system.

For the $n$-dimensional symmetry algebra $\mathcal{G}$, we firstly compute
the commutation relations between all the vector fields $v_i$ and $v_j$, which can be shown in a table, the entry in row $i$ and column $j$ representing $[v_i,v_j]=v_iv_j-v_jv_i$. Then taking any  subgroup $g=e^w (w=
\sum\limits_{j=1} ^{n} b_jv_j)$ to act on $v$, we have

\begin{eqnarray}
Ad_{exp(\epsilon w)}(v)&=& e^{-\epsilon w} v e^{\epsilon w} \nonumber \\
                       &=&v-\epsilon [w,v]+ \frac{1}{2!}\epsilon^2[w,[w,v]]-\cdots \label{ad1}\\
                       &=&(a_1v_1+\cdots+a_nv_n)-\epsilon [b_1v_1+\cdots+b_nv_n,a_1v_1+\cdots+a_nv_n]+o(\epsilon) \nonumber\\
                      &=&(a_1v_1+\cdots+a_nv_n)-\epsilon(\Theta_1v_1+\ldots+\Theta_nv_n)
                      +o(\epsilon), \nonumber
\end{eqnarray}
where $\Theta_i \equiv \Theta_i(a_1,\cdots, a_n, b_1,\cdots, b_n)$ can be easily obtained from the commutator table.

Equivalently, omitting $v_i$ we can rewrite (\ref{ad1}) as
$$(a_1,a_2,\cdots, a_n)  \xrightarrow{Ad_{exp(\epsilon w)}} (a_1-\epsilon\Theta_1,a_2-\epsilon\Theta_2,\cdots, a_n-\epsilon\Theta_n)+o(\epsilon).$$
According to the definition of the invariant, it is necessary that
\begin{equation}
\phi(a_1,a_2,\cdots, a_n)=\phi(a_1-\epsilon\Theta_1+o(\epsilon),a_2-\epsilon\Theta_2+o(\epsilon),\cdots, a_n-\epsilon\Theta_n+o(\epsilon)) \label{invariant}
\end{equation}
for any $b_i$.

Taking the derivative of Eq.~(\ref{invariant}) with respect to $\epsilon$ and
setting $\epsilon=0$, then extracting the coefficients of all $b_i$, $N (N \leq n)$ linear differential equations of $\phi$ are obtained. By solving these equations, all the invariants can be found.

\subsection{Calculation  of  the  adjoint  transformation  matrix}

The second task is the construction of the general adjoint  transformation matrix $A$, which is the product of the matrices of the separate adjoint actions $A_1, A_2, \cdots, A_n$. For further details, one can  refer to Ref.~\cite{Coggeshall} which showed three methods of constructing the adjoint matrix $A$. Here, before constructing the matrix $A$, one are able to draw a table, where the $(i,j)$-th entry gives $Ad_{exp(\epsilon v_i)}(v_j)$.

Firstly, applying the adjoint action of $v_1$  to $v=
\sum\limits_{i=1} ^{n} a_iv_i$ and with the help of adjoint representation table, we have
\begin{equation}
\label{nn}
\begin{aligned}
&Ad_{exp(\epsilon_1 v_1)}(a_1v_1+a_2v_2+\cdots+a_nv_n)\\
&=a_1Ad_{exp(\epsilon_1 v_1)}v_1+a_2Ad_{exp(\epsilon_1 v_1)}v_2+\cdots+a_nAd_{exp(\epsilon_1 v_1)}v_n\\
&=R_1v_1+R_2v_2+\cdots+R_nv_n,
\end{aligned}
\end{equation}
with $R_i\equiv R_i(a_1,a_2,\cdots, a_n,\epsilon_1), i=1\cdots n$.
To be intuitive,  the formula (\ref{nn}) can be rewritten into the following matrix form:
$$v \doteq (a_1,a_2,\cdots, a_n) \xrightarrow
{Ad_{exp(\epsilon_1 v_1)}} (R_1,R_2,\cdots, R_n)=(a_1,a_2,\cdots, a_n)A_1.
$$

Similarly, we can construct the  matrices  $A_2, A_3,\cdots A_n$ of  the  separate  adjoint actions of $v_2, v_3,\cdots, v_n$, respectively. Then the general adjoint  transformation  matrix $A$ is the  product  of  $A_1, \cdots, A_n$  taken  in  any  order
\begin{equation}
A\equiv A(\epsilon_1, \epsilon_2, \cdots, \epsilon_n)=A_1A_2\cdots A_n.
\end{equation}

That is to say, applying the most general adjoint action $Ad_{exp(\epsilon_n v_n)}\cdots Ad_{exp(\epsilon_2 v_2)}Ad_{exp(\epsilon_1 v_1)}
$ to $v$,  we have
\begin{equation}
\label{solveeq}
v \doteq (a_1,a_2,\cdots, a_n) \xrightarrow{Ad} (\tilde{a}_1,\tilde{a}_2,\cdots, \tilde{a}_n)=(a_1,a_2,\cdots, a_n)A.
\end{equation}

\emph{Remark 1:} In fact,  the right hand side of (\ref{solveeq}) can be regarded as $n$ algebraic
equations with respect to $\epsilon_1,\ldots, \epsilon_n$, which read
\begin{equation}
(\tilde{a}_1,\tilde{a}_2,\cdots, \tilde{a}_n)=(a_1,a_2,\cdots, a_n)A. \label{eqs}
\end{equation}
If Eqs.(\ref{eqs}) have the solution, it means that $\sum\limits_{i=1} ^{n} a_iv_i$ is equivalent to $\sum\limits_{i=1} ^{n} \tilde{a}_iv_i$ under the adjoint action; If Eqs.(\ref{eqs}) are incompatible, it shows that $\sum\limits_{i=1} ^{n} a_iv_i$  and $\sum\limits_{i=1} ^{n} \tilde{a}_iv_i$ are inequivalent.

\subsection{the classification of $\mathcal{G}$}

(1) \textbf{\emph{The first step:}}  scale the invariants.

If two vectors $v$ and $w$ are  adjoint equivalent, it is necessary that $\phi(v)=\phi(w)$ for any invariant $\phi$. However, if $v=cw$,  where $v$ and $w$ are also equivalent, their corresponding invariants satisfy $\phi(v)=c'\phi(w)$ and it is usually $\phi(v)\neq \phi(w)$. To avoid the latter case, we firstly make a scale to the invariant by adjusting  the coefficients of generators. Without loss of generality, one just need consider the values of the invariants to be $1$, $-1$ and $0$. To illustrate the point more clearly, we give three remarks.\\
\emph{Remark 2}: If the degree of the invariant is odd, we obtain $\phi(v)=c^{2k+1}\phi(w)$ with $v=cw$, then  the right $c$ can be selected to transform the  positive (negative) invariant into the  negative (positive) one. Now we just need consider two cases: $\phi=0$ and $\phi \neq0$
(for simplicity scaling it to 1 or -1).\\
\emph{Remark 3}: If the degree of the invariant is even (excluding zero), there is $\phi(v)=c^{2k}\phi(w)$ with $v=cw$, then  we can not choose the right $c$ to transform the  positive (negative) invariant into the  negative (positive) one. Now one  need consider three cases: $\phi=0$, $\phi>0$ and $\phi<0$.
Without loss of generality, we let $\phi=0$, $\phi=1$ and $\phi=-1$.\\
\emph{Remark 4}: Once one of the invariants is scaled (not zero), the other invariants (if any) can not be adjusted.

(2) \textbf{\emph{The second step:}} select the representative element.

According to different values of the invariants given in step 1, select the corresponding representative element in the simplest form named $\tilde{v}=\sum\limits_{i=1} ^{n} \tilde{a}_iv_i$.  Then solve  the adjoint transformation equation (\ref{eqs}).
If Eqs.(\ref{eqs}) have the solution with respect to $\epsilon_1, \cdots \epsilon_n$, it signifies that the selected representative element is right;
If Eqs.(\ref{eqs}) have no solution, we need reselect the proper representative element. Repeat the process until all the cases are finished in step 1.

\section{the new approach for the KdV equation and the heat equation}\label{heat equation}
\subsection{one-dimensional optimal system for the KdV equation}

The KdV equation  reads
\begin{equation}
\label{kdv}
u_t+u_{xxx}+uu_x=0,
\end{equation}
which arises in the theory of long waves in shallow water and other
physical systems in which both nonlinear and dispersive effects are relevant.
Using the classical Lie group method, one can obtain the symmetry algebra of (\ref{kdv}), i.e.
\begin{equation}
\label{kdvkdv}
v_1=\partial_{x},\quad  v_2=\partial_{t}, \quad
v_3=t\partial_x+\partial_u,\quad
v_4=x\partial_{x}+3t\partial_t-2u\partial_u.
\end{equation}
\emph{Step 1: calculate the invariants.}

The commutation relations between these vector fields is given by the following table, the entry in row $i$ and column $j$ representing $[v_i,v_j]=v_iv_j-v_jv_i$:
\begin{table}[htbp]
\centering
\caption{\label{table1}the commutator table of the KdV equation}
\begin{tabular}{c|cccc}
\hline
    &  $v_1$  &  $v_2$  &  $v_3$  &  $v_4$  \\
\hline
$v_1$ &  0    &  0    &  0     &  $v_1$ \\
\hline
$v_2$ &  0    &  0    &  $v_1$    &  3$v_2$ \\
\hline
$v_3$ &  0    &  $-v_1$    &  0    & $-2v_3$  \\
\hline
$v_4$ &  $-v_1$   &  $-3v_2$    &  $2v_3$    &  0 \\
\hline
\end{tabular}
\end{table}

Substituting $v=\sum\limits_{i=1} ^{4} a_iv_i$ and $w=\sum\limits_{j=1} ^{4} b_jv_j$ into (\ref{ad1}), there is
\begin{eqnarray}
Ad_{exp(\epsilon w)}(v) =(a_1v_4+\cdots+a_4v_4)-\epsilon(\Theta_1v_1+\ldots+\Theta_4v_4)
                      +o(\epsilon) \nonumber
\end{eqnarray}
with
\begin{equation}
\begin{aligned}
&\Theta_1=b_1a_4+b_2a_3-b_3a_2-b_4a_1,\quad
\Theta_2=3b_2a_4-3b_4a_2,\quad
\Theta_3=-2b_2a_4+2b_4a_3,\quad
\Theta_4=0.
\end{aligned}
\end{equation}
For any $b_i (i=1\cdots4)$, we have
\begin{equation}
\phi(a_1,a_2,a_3, a_4)=\phi(a_1-\epsilon\Theta_1+o(\epsilon),a_2-\epsilon\Theta_2+o(\epsilon),
a_3-\epsilon\Theta_3+o(\epsilon),a_4-\epsilon\Theta_4+o(\epsilon)). \label{invariant2}
\end{equation}

In Eq.~(\ref{invariant2}), taking the derivative of  with respect to $\epsilon$ and setting $\epsilon=0$, then extracting the coefficients of all $b_i$, four differential equations about $\phi(a_1,a_2,a_3,a_4)$ are directly obtained:
\begin{equation}
\label{kdvbubian}
  \left\{
   \begin{aligned}
&a_4\frac{\partial \phi}{\partial a_1}=0,\\
&a_3\frac{\partial \phi}{\partial a_1}+3a_4\frac{\partial \phi}{\partial a_2}=0,\\
&a_2\frac{\partial \phi}{\partial a_1}+2a_4\frac{\partial \phi}{\partial a_3}=0,\\
&a_1\frac{\partial \phi}{\partial a_1}+3a_2\frac{\partial \phi}{\partial a_2}-2a_3\frac{\partial \phi}{\partial a_3}=0.
\end{aligned}
\right.
\end{equation}

By solving Eqs.~(\ref{kdvbubian}), we obtain $\phi(a_1,a_2,a_3,a_4)=F(a_4)$, where $F$ is an arbitrary function of $a_4$. Here the basic invariant
of the KdV equation is only one, i.e. $a_4$, and $a_4$ is just the killing form given by Olver~\cite{Olver}.

\emph{Step 2: calculate the adjoint  matrix $A$.}

The adjoint representation table is given as
\begin{table}[htbp]
\centering
\caption{\label{table2}the adjoint representation table of the KdV equation
}
\begin{tabular}{c|cccc}
\hline
 Ad   &  $v_1$  &  $v_2$  &  $v_3$  &  $v_4$  \\
\hline
$v_1$ &  $v_1$    &  $v_2$    &  $v_3$    &  $v_4-\epsilon v_1 $ \\
\hline
$v_2$ &  $v_1$    &  $v_2$     & $v_3-\epsilon v_1$ &  $v_4-3\epsilon v_2$ \\
\hline
$v_3$ &  $v_1$    &  $v_2+\epsilon v_1$     &  $v_3$  &   $v_4+2\epsilon v_3$  \\
\hline
$v_4$ &  $e^{\epsilon}v_1$   &  $e^{3\epsilon}v_2$ & $e^{-2\epsilon}v_3$    &  $v_4$    \\
\hline
\end{tabular}
\end{table}

Applying the adjoint action of $v_1$ to
\begin{equation}
v=a_1v_1+a_2v_2+a_3v_3+a_4v_4, \label{kdvgvector}
\end{equation}
there is
\begin{equation}
\label{kdva1}
\begin{aligned}
Ad_{exp(\epsilon_1 v_1)}(a_1v_1+a_2v_2+a_3v_3+a_4v_4)
=(a_1-a_4\epsilon_1)v_1+a_2v_2+a_3v_3+a_4v_4
\end{aligned}
\end{equation}
It is easy to obtain
\begin{equation}
A_1=\left(
\begin{array}{cccccc}
 1 & 0 & 0 & 0  \\
 0 & 1 & 0 & 0\\
 0 & 0 & 1 & 0  \\
 -\epsilon_1 & 0 & 0 & 1
\end{array}
\right).
\end{equation}
Similarly,  we obtain $A_2, A_3$ and $A_4$:
\begin{equation}
A_2=\left(
\begin{array}{cccccc}
 1 & 0 & 0 & 0  \\
 0 & 1 & 0 & 0 \\
-\epsilon_2 & 0 & 1 & 0 \\
0 &-3\epsilon_2 & 0 & 1  \\
\end{array}
\right),\quad
A_3=\left(
\begin{array}{cccccc}
 1 & 0 & 0 & 0  \\
\epsilon_3 & 1 & 0 & 0 \\
0 & 0 & 1 & 0 \\
0 &0 &2\epsilon_3 & 1  \\
\end{array}
\right),\quad
A_4=\left(
\begin{array}{cccccc}
 e^{\epsilon_4} & 0 & 0 & 0  \\
 0 & e^{3\epsilon_4} & 0 & 0 \\
0 & 0 & e^{-2\epsilon_4} & 0 \\
0 &0 & 0 & 1  \\
\end{array}
\right).
\end{equation}
Then the general adjoint  transformation  matrix $A$ is obtained
\begin{equation}
A= A_1A_2A_3A_4=\left(
\begin{array}{cccc}
 e^{\epsilon_4} & 0 & 0 & 0  \\
\epsilon_3e^{\epsilon_4}& e^{3\epsilon_4} & 0 & 0 \\
-\epsilon_2e^{\epsilon_4} & 0 & e^{-2\epsilon_4} & 0 \\
(-\epsilon_1-3\epsilon_2\epsilon_3)e^{\epsilon_4}&-3\epsilon_2e^{3\epsilon_4} & 2\epsilon_3e^{-2\epsilon_4}  & 1  \\
\end{array}
\right).
\end{equation}

\emph{Step 3: the classification of symmetry algebra (\ref{kdvkdv}).}

According to ``\emph{Remark 2}", we have two cases: $a_4=1$ and $a_4=0$.  The adjoint transformation equations (\ref{eqs}) become
\begin{equation}
(\tilde{a}_1,\tilde{a}_2,\tilde{a}_3,\tilde{a}_4)=(a_1,a_2,a_3,a_4)A. \label{kdveqs}
\end{equation}
\emph{Case 1:} $a_4=1.$

Select a representative element $\tilde{v}=v_4$. Then by solving eqs.(\ref{kdveqs}),
we obtain the solution
\begin{equation}
\epsilon_1=a_1-\frac{1}{3}a_2a_3, \quad \epsilon_2=\frac{1}{3}a_2,\quad
\epsilon_3=-\frac{1}{2}a_3.
\end{equation}
That is to say, all the $v_4+a_1v_1+a_2v_2+a_3v_3$ are equivalent to $v_4$.\\
\emph{Case 2:} $a_4=0.$

Substituting $a_4=0$ into eqs.(\ref{kdvbubian}), we obtain a new invariant
$\phi(a_1,a_2,a_3)=a^2_2a^3_3$. In terms of ``\emph{Remark 2}", there are
also two cases: $a^2_2a^3_3=1$ and $a^2_2a^3_3=0$.

\emph{Case 2.1:} $a^2_2a^3_3=1$.

Adopt  two representative elements $\tilde{v}=v_2+v_3$ and $\tilde{v}=-v_2+v_3$.

For $a_2>0$ and $\tilde{v}=v_2+v_3$, eqs.(\ref{kdveqs}) with $a^2_2a^3_3=1$ have the solution
\begin{equation}
\epsilon_2=0, \quad \epsilon_3=-\frac{a_1}{a_2}, \quad \epsilon_4=-\frac{1}{3}\ln(a_2).
\end{equation}

For $a_2<0$ and $\tilde{v}=-v_2+v_3$, eqs.(\ref{kdveqs}) with $a^2_2a^3_3=1$  have the solution
\begin{equation}
\epsilon_2=0, \quad \epsilon_3=-\frac{a_1}{a_2}, \quad \epsilon_4=-\frac{1}{3}\ln(-a_2).
\end{equation}

\emph{Case 2.2:} $a^2_2a^3_3=0$.

(1) $a_3\neq0$ and $a_2=0$.

Adopt  two representative elements $v_3$ and $-v_3$. Then $\{a_1v_1+a_3v_3\}$ with $a_3>0$ is equivalent to $v_3$, where the solution for eqs.(\ref{kdveqs}) is
$\{\epsilon_2=\frac{a_1}{a_3}, \epsilon_4=\frac{1}{2}\ln(a_3)\}$, while $\{a_1v_1+a_3v_3\}$ with $a_3<0$ is equivalent to $-v_3$, where the solution for eqs.(\ref{kdveqs}) is
$\{\epsilon_2=\frac{a_1}{a_3}, \epsilon_4=\frac{1}{2}\ln(-a_3)\}$. Essentially,
$v_3$ and $-v_3$ is equivalent.

(2) $a_3=0$ : $a_2\neq0$ and $a_2=0$.

When $a_2\neq0$, similar to the above case (1), $a_2v_2+a_1v_1$ is equivalent to $v_2$
and $-v_2$.

When $a_2=0$, $a_1v_1$ is equivalent to $v_1$.

Recapitulating, an one-dimensional optimal system of  symmetry algebra (\ref{kdvkdv}) contains
\begin{equation}
\label{kdvop}
v_4; \quad v_3+v_2; \quad v_3-v_2; \quad v_3; \quad v_2; \quad v_1.
\end{equation}

The  optimal system obtained in (\ref{kdvop}) is just the same to that
found by Olver~\cite{Olver}.

\subsection{one-dimensional optimal system for the heat equation}
The equation for the conduction of heat in a one-dimensional road is written as
\begin{equation}
u_t=u_{xx}.                         \label{}
\end{equation}
The Lie algebra of infinitesimal symmetries for this equation is
spanned by  six vector fields
\begin{equation}
\label{vector}
\begin{aligned}
&v_1=\partial_{x},\quad  v_2=\partial_{t}, \quad
v_3=u\partial_u,\quad
v_4=x\partial_{x}+2t\partial_t,\\
&v_5=2t\partial_{x}-xu\partial_{u}, \quad
v_6=4tx\partial_{x}+4t^2\partial_{t}-(x^2+2t)u\partial_{u},
\end{aligned}
\end{equation}
and the infinitesimal subalgebra
\begin{equation}
v_{h}=h(x,t)\partial_{u},                    \nonumber
\end{equation}
where $h(x,t)$ is an arbitrary solution of the heat equation.
Since the infinite-dimensional subalgebra $\langle{v_h}\rangle$ does not lead to group invariant solutions, it will not be considered in the classification problem.

\emph{Step 1: calculate the invariants.}

Now consider the six-dimensional symmetry algebra $\mathcal{G}$ generated by $\{v_1,v_2, \cdots, v_6\}$ in (\ref{vector}).  Their commutator table is
given in table 3.
\begin{table}[htbp]
\centering
\caption{\label{table3}the commutator table of the heat equation}
\begin{tabular}{c|cccccc}
\hline
    &  $v_1$  &  $v_2$  &  $v_3$  &  $v_4$  & $v_5$    &  $v_6$ \\
\hline
$v_1$ &  0    &  0    &  0    &  $v_1$  &  $-v_3$  &  2$v_5$ \\
\hline
$v_2$ &  0    &  0    &  0    &  2$v_2$ &  2$v_1$  &  $4v_4-2v_3$ \\
\hline
$v_3$ &  0    &  0    &  0    &  0    &   0    &  0    \\
\hline
$v_4$ &  $-v_1$   &  $-2v_2$    &  0    &  0    &   $v_5$    &  2$v_6$    \\
\hline
$v_5$ &  $v_3$   &   $-2v_1$    &  0    &  $-v_5$    &   0    &  0    \\
\hline
$v_6$ &  $-2v_5$   &  $2v_3-4v_4$    &  0    &  $-2v_6$    &   0    &  0    \\
\hline
\end{tabular}
\end{table}

Substituting $v=\sum\limits_{i=1} ^{6} a_iv_i$ and $w=\sum\limits_{j=1} ^{6} b_jv_j$ into (\ref{ad1}), we have
\begin{eqnarray}
Ad_{exp(\epsilon w)}(v) =(a_1v_1+\cdots+a_6v_6)-\epsilon(\Theta_1v_1+\ldots+\Theta_6v_6)
                      +o(\epsilon) \nonumber
\end{eqnarray}
with
\begin{equation}
\begin{aligned}
&\Theta_1=-b_4a_1-2b_5a_2+b_1a_4+2b_2a_5,\quad
\Theta_2=-2b_4a_2+2b_2a_4,\quad
\Theta_3=b_5a_1+2b_6a_2-b_1a_5-2b_2a_6,\\
&\Theta_4=-4b_6a_2+4b_2a_6,\quad
\Theta_5=-2b_6a_1-b_5a_4+b_4a_5+2b_1a_6,\quad
\Theta_6=-2b_6a_4+2b_4a_6.
\end{aligned}
\end{equation}
Then consider the equation as follows
\begin{equation}
\phi(a_1,a_2,\cdots, a_6)=\phi(a_1-\epsilon\Theta_1+o(\epsilon),a_2-\epsilon\Theta_2+o(\epsilon),\cdots, a_6-\epsilon\Theta_6+o(\epsilon)) \label{invariant1}
\end{equation}
for any $b_i$.
Taking the derivative of Eq.~(\ref{invariant1}) with respect to $\epsilon$ and then setting $\epsilon=0$, extracting the coefficients of all $b_i$, five differential equations about $\phi(a_1,a_2,\cdots, a_6)$ are directly obtained:
\begin{equation}
\label{eqinvariant}
  \left\{
   \begin{aligned}
a_4 \frac{\partial \phi}{\partial a_1}-a_5 \frac{\partial \phi}{\partial
a_3}+2a_6 \frac{\partial \phi}{\partial
a_5}=0,  \\
a_4 \frac{\partial \phi}{\partial a_2}+a_5 \frac{\partial \phi}{\partial a_1}+a_6(2\frac{\partial \phi}{\partial
a_4}-\frac{\partial \phi}{\partial
a_3})=0,\\
-a_1\frac{\partial \phi}{\partial
a_1}-2a_2\frac{\partial \phi}{\partial
a_2}+a_5\frac{\partial \phi}{\partial
a_5}
+2a_6\frac{\partial \phi}{\partial
a_6}=0,\\
a_1\frac{\partial \phi}{\partial a_3}-2a_2\frac{\partial \phi}{\partial
a_1}-a_4\frac{\partial \phi}{\partial
a_5}=0,\\
-a_1\frac{\partial \phi}{\partial
a_5}+(2a_2\frac{\partial \phi}{\partial
a_3}-\frac{\partial \phi}{\partial
a_4})-
a_4\frac{\partial \phi}{\partial
a_6}=0.
   \end{aligned}
   \right.
  \end{equation}
Solving Eqs.~(\ref{eqinvariant}), one can obtain two basic common invariants
\begin{equation}
\Delta_1\equiv \phi_1(a_1,a_2,\cdots, a_6)=a_4^2-4a_2a_6
\end{equation}
and
\begin{equation}
\Delta_2\equiv \phi_2(a_1,a_2,\cdots,a_6)=a_4^3+2a_3a_4^2-4a_4a_2a_6+2a_4a_1a_5-8a_2a_3a_6-2a_2a_5^2-2a_1^2a_6.
\end{equation}
Here $\Delta_1$ is just the famous killing form in Ref.~\cite{Olver} while $\Delta_2$ is a new invariant of $v$ which is never addressed before.

\emph{Step 2: calculate the adjoint  matrix $A$.}

The adjoint representation table is given in table 4.
\begin{table}[htbp]
\centering
\caption{\label{table4}the adjoint representation
table of the heat equation}
\begin{tabular}{c|cccccc}
\hline
 Ad   &  $v_1$  &  $v_2$  &  $v_3$  &  $v_4$  & $v_5$    &  $v_6$ \\
\hline
$v_1$ &  $v_1$    &  $v_2$    &  $v_3$    &  $v_4-\epsilon v_1 $  &  $v_5+\epsilon v_3$  &  $v_6-2\epsilon v_5-\epsilon^2 v_3$ \\
\hline
$v_2$ &  $v_1$    &  $v_2$     &  $v_3$     & $v_4-2\epsilon v_2$ &  $v_5-2\epsilon v_1$  &  $v_6-4\epsilon v_4+2\epsilon v_3+4\epsilon^2v_2$ \\
\hline
$v_3$ &  $v_1$    &  $v_2$     &  $v_3$  &   $v_4$    &   $v_5$    &  $v_6$   \\
\hline
$v_4$ &  $e^{\epsilon}v_1$   &  $e^{2\epsilon}v_2$    &  $v_3$    &  $v_4$    &   $e^{-\epsilon}v_5$    &  $e^{-2\epsilon}v_6$    \\
\hline
$v_5$ &  $v_1-\epsilon v_3$   &   $v_2+2\epsilon v_1-\epsilon^2 v_3$    & $v_3$ & $v_4+\epsilon v_5$    &  $v_5$    &  $v_6$    \\
\hline
$v_6$ &  $v_1+2\epsilon v_5$   &  $v_2-2\epsilon v_3+4\epsilon
v_4+4\epsilon^2 v_6$  & $v_3$ &  $v_4+2\epsilon v_6$    &  $v_5$    &  $v_6$  \\
\hline
\end{tabular}
\end{table}
Applying the adjoint action of $v_1$ to
\begin{equation}
v=a_1v_1+a_2v_2+a_3v_3+a_4v_4+a_5v_5+a_6v_6, \label{gvector}
\end{equation}
there is
\begin{equation}
\label{a1}
\begin{aligned}
&Ad_{exp(\epsilon_1 v_1)}(a_1v_1+a_2v_2+a_3v_3+a_4v_4+a_5v_5+a_6v_6)\\
&=(a_1-a_4\epsilon_1)v_1+a_2v_2+(a_3+a_5\epsilon_1-a_6\epsilon_1^2)v_3+a_4v_4+(a_5-2\epsilon_1 a_6)v_5+a_6v_6,
\end{aligned}
\end{equation}
i.e.
$$v \doteq (a_1,a_2,\cdots, a_6) \xrightarrow{Ad_{exp(\epsilon_1 v_1)}} (a_1,a_2,\cdots, a_6)A_1.
$$
It is easy to obtain \begin{equation}
A_1=\left(
\begin{array}{cccccc}
 1 & 0 & 0 & 0 & 0 & 0 \\
 0 & 1 & 0 & 0 & 0 & 0\\
 0 & 0 & 1 & 0 & 0 & 0 \\
 -\epsilon_1 & 0 & 0 & 1 & 0 & 0 \\
 0 & 0 & \epsilon_1 & 0 & 1 & 0 \\
 0 & 0 & \epsilon_1^2 & 0 & -2\epsilon_1 & 0
 \end{array}
\right).
\end{equation}
Similarly, $A_2, \cdots, A_6$ are found to be
\begin{equation}
A_2=\left(
\begin{array}{cccccc}
 1 & 0 & 0 & 0 & 0 & 0 \\
 0 & 1 & 0 & 0 & 0 & 0\\
 0 & 0 & 1 & 0 & 0 & 0 \\
 0 &-2\epsilon_2 & 0 & 1 & 0 & 0 \\
 -2\epsilon_2 & 0 & 0 & 0 & 1 & 0 \\
 0 & 4\epsilon_2^2 & 2\epsilon_2 & -4\epsilon_2 & 0 & 1
 \end{array}
\right),\quad
A_4=\left(
\begin{array}{cccccc}
 e^{\epsilon_4} & 0 & 0 & 0 & 0 & 0 \\
 0 & e^{2\epsilon_4} & 0 & 0 & 0 & 0\\
 0 & 0 & 1 & 0 & 0 & 0 \\
 0 & 0 & 0 & 1 & 0 & 0 \\
 0 & 0 & 0 & 0 & e^{-\epsilon_4} & 0 \\
 0 & 0 & 0 & 0 & 0 & e^{-2\epsilon_4}
 \end{array}
\right).
\end{equation}

\begin{equation}
A_5=\left(
\begin{array}{cccccc}
 1 & 0 & -\epsilon_5 & 0 & 0 & 0 \\
 2\epsilon_5 & 1 & -\epsilon_5^2 & 0 & 0 & 0\\
 0 & 0 & 1 & 0 & 0 & 0 \\
 0 &0 & 0 & 1 & \epsilon_5 & 0 \\
 0 & 0 & 0 & 0 & 1 & 0 \\
 0 &0 & 0 & 0 & 0 & 1
 \end{array}
\right),\quad \quad  \quad \quad
A_6=\left(
\begin{array}{cccccc}
1  & 0 & 0 & 0 & 2\epsilon_6 & 0 \\
 0 & 1 & -2\epsilon_6 & 4\epsilon_6 & 0 & 4\epsilon_6^2\\
 0 & 0 & 1 & 0 & 0 & 0 \\
 0 & 0 & 0 & 1 & 0 & 2\epsilon_6 \\
 0 & 0 & 0 & 0 & 1 & 0 \\
 0 & 0 & 0 & 0 & 0 & 1
 \end{array}
\right),
\end{equation}
with $A_3=E$ being the identity matrix.

Hence the general adjoint  transformation  matrix $A$ is taken as
\begin{eqnarray}
A &=& A_4A_5A_3A_1A_2A_6\\
  &=& \left(
\begin{array}{cccccccc}
 e^{\epsilon_4}          & 0 & -\epsilon_5e^{\epsilon_4} & 0 & 2\epsilon_6 e^{\epsilon_4} & 0\\
 2\epsilon_5 e^{\epsilon_{4}} & e^{2\epsilon_4} & -(\epsilon_5^2+2\epsilon_6)e^{2\epsilon_4} & 4\epsilon_6e^{2\epsilon_4} & 4\epsilon_5\epsilon_6 e^{2\epsilon_4} & 4\epsilon_6^2 e^{2\epsilon_4}\\
 0 & 0 & 1 & 0 & 0 & 0 \\
 -\epsilon_1-2\epsilon_2\epsilon_5& -2\epsilon_2  & \epsilon_1\epsilon_5+4\epsilon_2\epsilon_6 & 1-8\epsilon_2\epsilon_6& -2\epsilon_1\epsilon_6- \epsilon_5\Xi& -2\epsilon_6\Xi \\
 -2\epsilon_2e^{-\epsilon_4} & 0 & \epsilon_1e^{-\epsilon_4} & 0 & -\Xi e^{-\epsilon_4}& 0\\
4\epsilon_1\epsilon_2e^{-2\epsilon_4}&
4\epsilon_2^2e^{-2\epsilon_4}& -(\epsilon_1^2+2\epsilon_2\Xi)e^{-2\epsilon_4}
&
4\epsilon_2\Xi e^{-2\epsilon_4}&
2\epsilon_1\Xi e^{-2\epsilon_4}&
\Xi^2
e^{-2\epsilon_4}
\end{array}
\right).
\end{eqnarray}
with $\Xi=4\epsilon_2\epsilon_6-1$.
%

\emph{Step 3: the classification of symmetry algebra (\ref{vector}).}

In the following, two invariants $\Delta_1$ and $\Delta_2$ will be made full use of to give an classification  of the  algebra $\mathcal{G}$.  Since the degree of $\Delta_1=a_4^2-4a_2a_6$  is two, we can
scale it to three cases: $\Delta_1=1, \Delta_1=-1$ and $\Delta_1=0$ according to ``\emph{Remark 3}".\\
\emph{Case 1:} $\Delta_1=a_4^2-4a_2a_6=1,  \Delta_2=c.$

Here $c$ is an arbitrary real constant. Under $\Delta_1=1$ and $\Delta_2=c$, choose a representative element, for example, select $\tilde{v}=v_4+\frac{c-1}{2}v_3$ (i.e. $\tilde{a}_4=1, \tilde{a}_3=\frac{c-1}{2}$).

From $\Delta_1=a_4^2-4a_2a_6=1$, we know that $a_2, a_4, a_6$ can not be all zeros simultaneously. Without loss of generality, one  only considers  $a_6 \neq 0$. For $a_6=0$ ($a_2 \neq 0$ or $ a_4 \neq 0$), one can transform it into the
the case of  $a_6 \neq 0$ by selecting the appropriate $\epsilon_i (i=1\cdots 6)$ which are shown in eqs.(\ref{eqs}).

For $a_6 \neq 0$, the general solution of $\Delta_1=1$ and $\Delta_2=c$
is
\begin{equation}
a_2=\frac{a_4^2-1}{4a_6}, \quad a_3=a_1^2a_6-a_1a_4a_5-\frac{1}{2}a_4+\frac{c}{2}-\frac{a_4^2-1}{4a_6}a_5^2. \label{case1a}
\end{equation}
where $a_1,a_4,a_5,a_6$ are arbitrary real constants. According to the
formula (\ref{eqs}), six algebra equations about $\epsilon_i$ are proposed. After
substituting $\tilde{a}_4=1, \tilde{a}_3=\frac{c-1}{2}, \tilde{a}_1=\tilde{a}_2=\tilde{a}_3=\tilde{a}_4=0$ with (\ref{case1a}) into these equations, one can find the solution:
\begin{equation}
\epsilon_1=\frac{a_5+2a_1a_4a_6-a_4^2a_5}{2a_6}e^{\epsilon_4},\quad
\epsilon_2=\frac{a_4-1}{4a_6}e^{2\epsilon_4},\quad
\epsilon_5=(2a_1a_6-a_4a_5)e^{-\epsilon_4},\quad
\epsilon_6=-\frac{1}{2}a_6e^{-2\epsilon_4}.
\end{equation}
\emph{Case 2:} $\Delta_1=a_4^2-4a_2a_6\equiv -1,  \Delta_2=c.$

From $\Delta_1=-1$, it illustrates $a_6\neq 0$. Now the relation among $a_i$ reads
\begin{equation}
a_2=\frac{a_4^2+1}{4a_6}, \quad a_3=-a_1^2a_6+a_1a_4a_5-\frac{1}{2}a_4-\frac{c}{2}-\frac{a_4^2+1}{4a_6}a_5^2. \label{case2}
\end{equation}

When $a_6>0$ and $a_6<0$, take the representative element $\tilde{v}=\frac{1}{2}(v_2+v_6-cv_3)$ and $\tilde{v}=\frac{1}{2}(-v_2-v_6-cv_3)$ respectively. Then Eqs.(\ref{eqs}) with (\ref{case2}) are separately proved right by selecting
\begin{equation}
\epsilon_1=-\frac{\sqrt{2}}{2}\frac{2a_1a_4a_6-a_5-a_4^2a_5}{\sqrt{a_6}},\quad
\epsilon_2=\frac{1}{2}a_4,\quad
\epsilon_4=\frac{1}{2}\ln(2a_6), \quad
\epsilon_5=-\frac{\sqrt{2}}{2}\frac{2a_1a_6-a_4a_5}{\sqrt{a_6}},\quad
 \epsilon_6=0
\end{equation}
and
\begin{equation}
\epsilon_1=-\frac{1}{2}\frac{\sqrt{-2a_6}(2a_1a_4a_6-a_5-a_4^2a_5)}{a_6},\quad
\epsilon_2=-\frac{1}{2}a_4,\quad
\epsilon_4=\frac{1}{2}\ln(-2a_6), \quad
\epsilon_5=\frac{1}{2}\frac{\sqrt{-2a_6}(2a_1a_6-a_4a_5)}{a_6},\quad
 \epsilon_6=0.
\end{equation}
In this case, the general one-dimensional Lie algebra (\ref{gvector}) is equivalent
to $v_2+v_6+\beta v_3$ with $\beta$ being arbitrary.\\
\emph{Case 3:} $\Delta_1=0,  \Delta_2=c.$

Notice that  $\Delta_1=0$  and $\Delta_2$ itself is an
odd polynomial with respect to $a_i$, so one just need to
consider $\Delta_2=1$ and $\Delta_2=0$ according to ``\emph{Remark 2}".

\emph{Case 3.1:} $\Delta_1=0,  \Delta_2=1.$

Select a representative element $\tilde{v}=-v_2-\frac{\sqrt{2}}{2}v_5$.

When $a_6\neq 0$, there must be $a_6<0$ for  the identity $2a_6=-(2a_1a_6-a_4a_5)^2$ solved by $\Delta_1=0$ and $\Delta_2=1$. Under the  restriction of invariants, we have
\begin{equation}
a_2=\frac{a_4^2}{4a_6},\quad a_1=\frac{a_4a_5+\sqrt{-2a_6}}{2a_6} \quad {\rm{and}} \quad   a_2=\frac{a_4^2}{4a_6},\quad a_1=\frac{a_4a_5-\sqrt{-2a_6}}{2a_6}.                           \label{case31}
\end{equation}
Then after choosing
\begin{equation}
\begin{aligned}
&\epsilon_1=\frac{\sqrt{2}}{4}\cdot\frac{e^{\epsilon_4}}{a_6}(\sqrt{2}a_5+e^{\epsilon_4}
+\sqrt{2}a_4\epsilon_5e^{\epsilon_4}),\quad
\epsilon_2=\frac{e^{\epsilon_4}}{4a_6}(\mp2\sqrt{-a_6}+a_4e^{\epsilon_4}),\\
&\epsilon_5=\pm \frac{\sqrt{2}}{8}\cdot\frac{e^{-\epsilon_4}}{\sqrt{-a_6}}(4a_4a_6+
2a_5^2+8a_3a_6+e^{2\epsilon_4}),\quad
\epsilon_6=\pm \frac{1}{2}\sqrt{-a_6}e^{-\epsilon_4},
\end{aligned}
\end{equation}
one can transform (\ref{gvector}) with (\ref{case31}) into $\tilde{v}=-v_2-\frac{\sqrt{2}}{2}v_5$.

When $a_6=0$, one can act by any $Ad_{exp(\epsilon_6 v_6)}$ with $\epsilon_6\neq0$ to get a nonzero coefficient in front of $v_6$, reducing to the previous case.

\emph{Case 3.2:} $\Delta_1=0,  \Delta_2=0.$

\emph{Case 3.2.1:} Not all $a_2, a_4$ and $a_6$ are zeros. Without loss of generality, we just consider the case of $a_6\neq 0.$

Substituting $\Delta_1=0$ and $\Delta_2=0$ into Eqs.(\ref{eqinvariant}), we  obtain a new invariant
\begin{equation}
\Delta_3=4a_3+2a_4+\frac{a_5^2}{a_6}.  \label{}
\end{equation}
Then there are now two cases, depending on the sign of the invariant $\Delta_3$:

(1) $\Delta_3=1$. Solving $\Delta_3=1$, we have
\begin{equation}
a_1=\frac{a_4a_5}{2a_6}, \quad a_2=\frac{a_4^2}{4a_6}, \quad a_3=\frac{a_6-a_5^2-2a_4a_6}{4a_6}.     \label{case3211}
\end{equation}

When $a_6>0$, choose the representative element $\tilde{v}=\frac{1}{4}v_3+v_6$,
then Eqs.(\ref{eqs}) with (\ref{case3211}) have the solution
\begin{equation}
 \epsilon_1=\frac{a_5\sqrt{a_6}+a_4a_6\epsilon_5}{2a_6},\quad \epsilon_2=\frac{a_4}{4},
\quad \epsilon_4=\frac{1}{2}\ln{a_6}.
\label{a61}
\end{equation}

When $a_6<0$, the representative element is taken as $\tilde{v}=\frac{1}{4}v_3-v_6$.
It is easy to see that Eqs.(\ref{eqs}) are right with
\begin{equation}
 \epsilon_1=\frac{a_5\sqrt{-a_6}-a_4a_6\epsilon_5}{2a_6},\quad \epsilon_2=-\frac{a_4}{4},
\quad \epsilon_4=\frac{1}{2}\ln{(-a_6)}.
\label{a62}
\end{equation}

Here it is noted that $\frac{1}{4}v_3+v_6$ and $\frac{1}{4}v_3-v_6$ is inequivalent.

(2) $\Delta_3=0$. Now we have
\begin{equation}
a_1=\frac{a_4a_5}{2a_6}, \quad a_2=\frac{a_4^2}{4a_6}. \label{}
\end{equation}

It can be easily proved that via the same adjoint transformation (\ref{a61}) and (\ref{a62}), the  Lie algebra (\ref{gvector}) is converted into
$v_6$ and $-v_6$, respectively.

\emph{Case 3.2.2:} $a_2=a_4=a_6=0$. Substituting $a_2=a_4=a_6=0$ into (\ref{eqs}), we find that it can also be divided into two cases:

(1) Not all $a_1$ and $a_5$ are zeros. Here we suppose $a_5\neq0$.

%
%

When $a_5\neq0$, give a representative element $v_1$. One can see that Eqs.(\ref{eqs}) with $a_2=a_4=a_6=0$  has a solution
\begin{equation}
\epsilon_1=\frac{e^{\epsilon_4}(e^{\epsilon_4}\epsilon_5a_1-a_3)}{a_5},\quad
\epsilon_2=\frac{1}{2}\frac{e^{\epsilon_4}(e^{\epsilon_4}a_1-1)}{a_5},\quad
\epsilon_6=-\frac{1}{2}a_5e^{-\epsilon_4}.
\label{}
\end{equation}

(2) $a_1=a_5=0$. Now we have $a_1=a_2=a_4=a_5=a_6=0$ and the  general Lie algebra (\ref{gvector}) becomes $v_3$.

In summary, an optimal system of one-dimensional subalgebras of the heat equation is found to be those spanned by
\begin{equation}
\label{heatop}
\begin{aligned}
&\omega_1(\alpha)=v_4+\alpha v_3,\quad \quad    \alpha \in \mathbb{R},\\
&\omega_2(\beta)=v_2+v_6+\beta v_3, \quad \quad   \beta \in \mathbb{R},\\
&\omega_3=v_2+\frac{\sqrt{2}}{2}v_5,\\
&\omega_4=\frac{1}{4}v_3+v_6,\\
&\omega_5=\frac{1}{4}v_3-v_6,\\
&\omega_6=v_6,\\
&\omega_7=v_1,\\
&\omega_8=v_3.
\end{aligned}
\end{equation}

The resulting  optimal system (\ref{heatop}) of the heat equation is really
optimal and completely equivalent to that given in Ref.~\cite{Qu1}, which  is a further reduction to the result of Olver~\cite{Olver}.\\

\textbf{\emph{Remark 5:}} The key point of our new method is to solve some algebra equations which are reflected in (\ref{eqs}) and it can easily be carried out by Maple.

%
%
%

\section{Summary and discussion}\label{conclusion}

Group invariant solutions have been used to great effect in the
description of the asymptotic behaviour of much more general
solutions to systems of partial differential equations. These group invariant solutions are characterized by their invariance under some symmetry group of the
system of partial differential equations. Since there are almost always an infinite
number of different symmetry groups one might employ to find group invariant solutions, a means of determining which groups give fundamentally different types of invariant solutions is essential for gaining a complete understanding of the solutions which might be available. This classification problem can be solved by looking at the adjoint representation of the symmetry group on its Lie algebra, which firstly used by Ovsiannikov. The one-dimensional classification of the symmetry algebras of the KdV equation and the heat equation are demonstrated by Olver with an easy-to-operate method in detail, which only depends on the fragments of the theory of Lie algebras. However, as Olver said, in essence this problem is attacked by the na\"{i}ve approach of taking a general element in Lie algebra and subjecting it to various adjoint transformations so as to ``simplify" it as much as possible. To make up this  problem and ensure the comprehensiveness with inequivalence,  we develop a direct and systemic algorithm for  the one-dimensional optimal system.  The new approach is very natural and every elements in the optimal system can be found step by step.

Our method  introduced in this paper, which is essentially new, only depends on the commutator and adjoint representative relations  among  the generators of Lie algebras. The main work includes:

(1) A valid method is proposed to compute all the general invariants of the one-dimensional  Lie algebra, which include the well-known  killing form;

(2) A criterion is introduced to scale the invariants, which appears in ``\emph{Remark 2}'', ``\emph{Remark 3}'' and ``\emph{Remark 4}'';

(3) For two one-dimensional subalgebras $v$ and $\tilde{v}$, we introduce
an algebraic equations system  (\ref{eqs})  to determine their equivalences in the sense of adjoint transformation;

(4) Based on all the scaled  invariants, we put forward
a direct and effective algorithm to construct one-dimensional optimal system. With the new approach, every  element in the optimal system can be  found step by step.

Since all the representative elements are attached to different values of the invariants, it ensures  the optimality of the optimal system. From the process of the operation in our method, one can easily see  that how  these representatives are mutually inequivalent. Due to the optimal system of symmetry algebra, a family of group invariant solutions can be recovered.
Since many important equations arising from physics are of low dimensions and reducing them to ODEs requires only the determination of small parameter optimal systems, we hope this method will be useful elsewhere. Furthermore, how to apply all the invariants to construct r-parameter ($r\geq 2$) optimal systems is in our  consideration. Since the algorithm is very systemic,
we believe that it will provide a very good manner for the mechanization.

\section*{Acknowledgments}
The authors extend their gratitude to
Professors Qu C Z and Lou S Y for their helpful discussions.
This work is supported by Zhejiang Provincial Natural Science Foundation of China under Grant No. LQ13A010014, the National Natural Science
Foundation of China (Grant Nos. 11326164 and 11275072), the Research Fund for the Doctoral Program of Higher Education of China (Grant No. 20120076110024) and the Innovative Research Team Program of the National Natural Science Foundation of China (Grant No. 61321064) and K C Wang Magna Fund in Ningbo University.

\bibliography{<your-bib-database>}

\end{document}